\documentclass[11pt]{article}

 \usepackage{epsfig}
\usepackage{algorithm}
\usepackage{algorithmic}
\usepackage{latexsym}
\usepackage{amsmath}
\usepackage{cite}
\usepackage{verbatim}
\usepackage{amsfonts}
\usepackage{graphicx}
\usepackage[normalem]{ulem}
\usepackage{array}
\usepackage{enumerate}
\usepackage{authblk}
\usepackage{color}\usepackage[dvipsnames]{xcolor}

\usepackage{tikz}\usetikzlibrary{patterns}
\usepackage{makecell}

\DeclareFontFamily{OT1}{pzc}{}
\DeclareFontShape{OT1}{pzc}{m}{it}{<-> s * [1.10] pzcmi7t}{}
\DeclareMathAlphabet{\mathpzc}{OT1}{pzc}{m}{it}

\usepackage[
  bookmarks=false, 
  colorlinks,
  citecolor=brown!70!black,
  linkcolor=brown!80!black,
  urlcolor=blue!70!black,
]{hyperref}
\usepackage{relsize,exscale}
\newtheorem{lem}{Lemma}

\newtheorem{cor}{Corollary}
\newtheorem{thm}{Theorem}
\newtheorem{rem}{Remark}

\title{Bijections from  Dyck and Motzkin meanders with catastrophes to pattern avoiding Dyck paths}

\author[]{Jean-Luc Baril\thanks{Corresponding author: \texttt{barjl@u-bourgogne.fr}} }
\author[]{Sergey Kirgizov}

\affil[]{LIB, Univ. Bourgogne Franche-Comt{\'e}, France}

\begin{document}
\maketitle

\begin{abstract} In this note, we present constructive bijections from  Dyck and Motzkin meanders with catastrophes to Dyck paths avoiding some patterns. As a byproduct, we deduce correspondences from Dyck and Motzkin excursions to restricted Dyck paths.
\end{abstract}
\section{Introduction and notations}
The domain of lattice paths provides a very fertile ground for the combinatorial community. They have many applications in computer science, queuing theory, biology and physics \cite{Sta}, and  there are a multitude of one-to-one correspondences  with various combinatorial objects such as directed animals, pattern avoiding  permutations,  bargraphs,  RNA structures  and so on \cite{Bar1,Knu,Sta}. A recurring problem in combinatorics is the enumeration of these paths with respect to their length and  other  statistics \cite{Barc,Bar,Deu,Man1,Mer,Pan,Sap,Sun}. In the literature, Dyck and Motzkin paths are the most often considered,  possibly because  they  are, respectively,  counted  by the  famous  Catalan and  Motzkin  numbers  (see \href{https://oeis.org/A108}{A108} and \href{https://oeis.org/A1006}{A1006} in the  Sloane's
On-line
Encyclopedia  of  Integer Sequences \cite{Sloa}).

Throughout this note, a {\it lattice path} is defined by a starting point $(0,0)$, an ending point $(n,k)$ with $n,k\geq 0$, it consists of steps lying in $S=\{(1,i): i\in\mathbb{Z},  i\leq 1\}$, and it never goes below the $x$-axis. The length of  a  path  is  the  number  of  its  steps. We denote  by $\epsilon$ the  empty path, {\it i.e.}, the path of length zero. Constraining the steps to be in  $\{(1,1),(1,-1)\}$ or $\{(1,1),(1,0),(1,-1)\}$,  and fixing the end point on the $x$-axis,  we retrieve the well-known  definition  of {\it Dyck  and Motzkin paths} \cite{Sta} respectively.  Let $\mathcal{D}_{n}$ be the set of Dyck paths  of semilength $n$, we define $\mathcal{D}=\cup_{n\geq 0}\mathcal{D}_{n}$. For short, we set $U=(1,1)$, $D=(1,-1)$, $F=(1,0)$ and
$D_i=(1,-i)$ for $i\geq 2$.

\begin{figure}[h]
\begin{center}($a$)
\scalebox{0.55}{\begin{tikzpicture}[ultra thick]
   \draw[black, line width=2pt]
(0,0)--(0.4,0.4)--(0.8,0.8)--(1.2,0.4)--(1.6,0.4)--(2,0.8)--(2.4,1.2)--(2.8,1.2)--(3.2,0)--(3.6,0.4)--(4,0)--(5.2,1.2)--(6,0.4)--(6.4,0.8)--(6.8,0)--(7.6,0.8)--(8,0.8)--(8.4,0.8)--(8.8,1.2)--(9.2,1.2);
   \draw  (3.3,1) node {\Large $D_3$};\draw  (6.8,1) node {\Large $D_2$};

   \tikzset{every node/.style={circle, white, fill = black, inner sep = 1.5pt}}
   \node at (0,0) {};
   \node at (0.4,0.4) {};
   \node at (0.8,0.8) {};
   \node at (1.2,0.4) {};
   \node at (1.6,0.4) {};
   \node at (1.6,0.4) {};
   \node at (2,  0.8) {};
   \node at (2.4, 1.2) {};
   \node at (2.8, 1.2) {};
   \node at (3.2, 0) {};
   \node at (3.6, 0.4) {};
   \node at (4, 0) {};
   \node at (4, 0) {};
   \node at (4.4, 0.4) {};
   \node at (4.8, 0.8) {};
   \node at (5.2, 1.2) {};
   \node at (5.6, 0.8) {};
   \node at (6, 0.4) {};
   \node at (6.4, 0.8) {};
   \node at (6.8, 0) {};
   \node at (7.2, 0.4) {};
   \node at (7.6, 0.8) {};
   \node at (8, 0.8) {};
   \node at (8.4, 0.8) {};
   \node at (8.8, 1.2) {};
   \node at (9.2, 1.2) {};
   \draw[black, thick] (0,0)--(9.2,0); \draw[black, thick] (0,0)--(0,1.8);\end{tikzpicture}}
  \qquad ($b$)
\scalebox{0.55}{\begin{tikzpicture}[ultra thick]
   \draw[black, line width=2pt]
   (0,0)--(0.4,0.4)--(0.8,0.8)--(1.2,0.4)--(2,1.2)--(2.4,0.8)--(2.8,1.2)--(3.2,0)--(3.6,0.4)--(4,0)--(5.2,1.2)--(6,0.4)--(6.4,0.8)--(6.8,0);
   \draw  (3.3,1) node {\Large $D_3$};\draw  (6.8,1) node {\Large $D_2$};

   \tikzset{every node/.style={circle, white, fill = black, inner sep = 1.5pt}}
   \node at (0,0) {};
   \node at (0.4,0.4) {};
   \node at (0.8,0.8) {};
   \node at (1.2,0.4) {};
   \node at (1.6,0.8) {};
   \node at (2,1.2) {};
   \node at (2.4,0.8) {};
   \node at (2.8,1.2) {};
   \node at (3.2,0) {};
   \node at (3.6,0.4) {};
   \node at (4,0) {};
   \node at (4.4,0.4) {};
   \node at (4.8,0.8) {};
   \node at (5.2,1.2) {};
   \node at (5.6,0.8) {};
   \node at (6,0.4) {};
   \node at (6.4,0.8) {};
   \node at (6.8,0) {};
   \draw[black, thick] (0,0)--(6.8,0);\draw[black, thick] (0,0)--(0,1.8);
\end{tikzpicture}}
\end{center}
\caption{($a$)  A Motzkin meander  with catastrophes in $\mathcal{M}'_{23}$, and ($b$) a Dyck  excursion with catastrophes in $\mathcal{E}_{17}$ }
\label{fig1}
\end{figure}
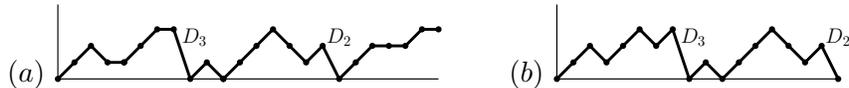

Considering these notations, a {\it Motzkin meander with catastrophes} is a lattice path where possible steps are $U, D, F$ and $D_i$ for $i\geq 2$, such that all steps $D_i$ end on the $x$-axis,  and if we add the property that the path ends on the $x$-axis, we call it a {\it Motzkin  excursion} with catastrophes  (see \cite{Ban}). {\it Dyck meanders} and {\it Dyck excursions} with catastrophes are those avoiding the step $F$.
Let $\mathcal{M}_n$ (resp. $\mathcal{E}_n$) be the set of length $n$ Dyck meanders (resp. excursions) with catastrophes, and we set  $\mathcal{M}=\cup_{n\geq 0}\mathcal{M}_n$ (resp. $\mathcal{E}=\cup_{n\geq 0}\mathcal{E}_n$).
The sets of Motzkin meanders and excursions with catastrophes are
respectively
denoted  by adding  prime superscripts, $\mathcal{M}'$ and $\mathcal{E}'$.  As mentioned in Corollary 2.4 in \cite{Ban}, the cardinality of $\mathcal{M}_n$ is given by the sequence \href{https://oeis.org/A274115}{A274115}  in \cite{Sloa}, and the cardinality of $\mathcal{E}_n$ is given by the sequence \href{https://oeis.org/A224747}{A224747}. For instance, we have $UUDFUUFD_3UDUUUDDUD_2UUFFUF\in \mathcal{M}'_{23}$ and $UUDUUDUD_3UDUUUDDUD_2\in \mathcal{E}_{17}$, and we refer to Figure \ref{fig1} for an illustration of these two paths.
Since Motzkin meanders with catastrophes can be obtained from Dyck meanders with catastrophes by possibly adding flat steps $F$, the ordinary generating function (o.g.f.) for the  cardinality of $\mathcal{M}'_n$  is given by $\frac{M(x/(1-x))}{1-x}$ where 
$M(x) = \frac{2 x}{2 x + \left(x + 1\right) \left(\sqrt{1 - 4 x^{2}} - 1\right)}$
is the o.g.f. for $\mathcal{M}_n$ (see \cite{Ban}), which generates the  $(n+1)$th term of \href{https://oeis.org/A54391}{A54391}. Simarly,  $\mathcal{E}'_n$ is   counted by the $n$th term of \href{https://oeis.org/A54391}{A54391}.

Dyck meanders with catastrophes was first introduced by Krinik {\em et al.}
in~\cite{Kri} in the context of queuing theory. They  correspond to the evolution of the queue by allowing some resets modeled by a catastrophe step $D_i$ for $i\geq 2$. 
Recently in \cite{Ban},  Banderier and Wallner provide many results about the enumeration and limit laws of these objects. Using  algrebraic methods  they prove that the set $\mathcal{M}_n$ of length $n$ Dyck meanders with catastrophes  has the same cardinality as the set of equivalence classes of semilength $n+1$ Dyck paths modulo the positions of the pattern $DUU$,  which in turn (see \cite{Mane}) is in one-to-one correspondence with the set $\mathcal{A}_{n}$ of semilength $n$ Dyck paths avoiding occurrences at height $h>0$ of the patterns $UUU$ and $DUD$. They also provide a constructive bijection between $\mathcal{E}_n$ and the set of length $n$ Motzkin paths having their flat steps $F$ at height one.

The motivation of this work is to exhibit one-to-one correspondences between restricted Dyck paths (with no catastrophes) and the sets of paths with catastrophes  $\mathcal{M}_n$, $\mathcal{E}_n$, $\mathcal{M}'_n$, and $\mathcal{E}'_n$. In Section 2 we present a constructive bijection between $\mathcal{M}_n$ and $\mathcal{A}_n$. Considering its restriction to  excursions with catastrophes, we prove that $\mathcal{E}_n$ is in one-to-one correspondence with the set $\mathcal{A}'_n$ of Dyck paths in $\mathcal{A}_n$ where  any occurrence $UD$ on the $x$-axis appears before an occurrence  of $UUU$  (not necessarily contiguous to the occurrence $UD$). This bijection establishes a curious correspondence connecting Dyck meanders with catastrophes and equivalence classes modulo the  positions of  $DUU$ in Dyck paths. In Section 3 we conduct the  counterpart study for Motzkin meanders and excursions. More precisely, we exhibit a bijection between $\mathcal{M}'_n$ and the set $\mathcal{B}_{n+1}$ of semilength $n+1$ Dyck paths avoiding the pattern $UUU$ at height $h\geq 2$, which also induces a bijection  from the set  $\mathcal{E}'_n$ of Motzkin excursions with catastrophes to the set $\mathcal{B}_{n}$.  The following table gives an overview of all these correspondences.

\begin{center}
\scalebox{0.96}{  \small
  \setlength{\tabcolsep}{5pt}
  \begin{tabular}{  r | c  c  c | l  }
    Dyck  meanders  with cat. &$\mathcal{M}_n$& $\to$ &  $\mathcal{A}_n$ 
    & 
    \makecell[cl]{
    Dyck paths avoiding $UUU$ \\
    and $DUD$ at $h>0$
    }
    \\
    \hline
    Dyck  excursions  with  cat. & $\mathcal{E}_n$  & $\to$ &  $\mathcal{A'}_n$ 
    &
    \makecell[cl]{
    $\mathcal{A}_n$ whose every $UD$ on the \\
    $x$-axis appears before $UUU$ 
    }
    \\ 
     \hline
    Dyck  paths & $\mathcal{D}_n$  & $\to$ &  $\mathcal{A}^\star_{2n}$ 
    &
    \makecell[cl]{
    Dyck paths starting with $UU$ \\
    and avoiding $UUU$ and $DUD$
    }
    \\\hline
    Motzkin meanders with cat. & $\mathcal{M}'_n$   & $\to$ &  $\mathcal{B}_{n+1}$
    & 
    \makecell[cl]{
    Dyck paths avoiding $UUU$ \\
    at $h\ge 2$
    }
    \\ \hline
    Motzkin excursions with cat. & $\mathcal{E}'_n$    & $\to$ &  $\mathcal{B}_{n}$
    &
    \\ \hline
    Motzkin  paths & 
    $\mathcal{M}\mathpzc{otz}_n $ 
    & $\to$ & 
    $\mathcal{B}'_{n+1}$
    & 
    \makecell[cl]{
    Dyck paths avoiding $UUU$ \\
    at $h \ge 2$ and $DU$ at $ h = 1$.
    }
    \\
    \end{tabular}}
\end{center}

\section{Dyck meanders with catastrophes}

In this section we exhibit a constructive bijection between the set $\mathcal{M}_n$ of length $n$ Dyck meanders with catastrophes and the set
 $\mathcal{A}_n$ of semilength $n$ Dyck paths having no occurrence of the consecutive three steps $UUU$ and $DUD$ at height $h>0$ (or equivalently with a minimal ordinate $h>0$).
We set $\mathcal{A}=\cup_{n\geq 0}\mathcal{A}_n$. Let us define recursively the map $\phi$ from $\mathcal{M}$ to  $\mathcal{D}$ as follows.
 For $P\in\mathcal{M}$, we set
 $$\phi(P)=\left\{\begin{array}{ll}
\epsilon& \mbox{ if } P=\epsilon,\qquad\hfill (i)\\
UD\phi(\alpha) &\mbox{ if } P=U\alpha,\qquad\hfill (ii)\\
UUD\phi(\alpha )D\phi(\beta)  & \mbox{ if } P= U\alpha D \beta, \qquad\hfill  (iii)\\
U\phi(\alpha D) D\phi(\beta)  & \mbox{ if } P= U\alpha D_2 \beta,\qquad \hfill (iv)\\
UD\phi(\alpha D_{i-1}) \phi(\beta)  & \mbox{ if } P= U\alpha D_i \beta \mbox{ and } i\geq 3,\qquad\hfill (v)\\
\end{array}\right.$$
where $\beta\in\mathcal{M}$, and $\alpha$ is either the empty path or a lattice path consisting of $U$- and $D$-steps such that $\alpha$ (resp. $\alpha D$, $\alpha D_{k-1}$) ends on the $x$-axis in the case ($iii$) (resp. ($iv$), ($v$)), and $\alpha$ does not necessarily end on the $x$-axis in the case ($ii$).

\begin{figure}[h]
\begin{center}
\scalebox{0.55}{\begin{tikzpicture}[ultra thick]
 \draw[black, thick] (0,0)--(7,0); \draw[black, thick] (0,0)--(0,2);
 \draw[black, dashed, very thick] (0.4,0.4)--(3,0.4);
  \draw[black, line width=3pt] (0,0)--(0.4,0.4);
 \draw[orange,very thick] (0.4,0.4) parabola bend (1.7,2) (3,0.4);
 \draw  (1.7,1) node {\huge $\mathbf{\alpha}$};
 \end{tikzpicture}}
 \qquad$\stackrel{\phi}{\longrightarrow}$\qquad
\scalebox{0.55}{\begin{tikzpicture}[ultra thick]
 \draw[black, thick] (0,0)--(7,0); \draw[black, thick] (0,0)--(0,2);
 \draw[black, line width=3pt] (0,0)--(0.4,0.4)--(0.8,0);
 \draw[orange,very thick] (0.8,0) parabola bend (1.9,1.5) (3,0);
 \draw  (1.9,0.7) node {\huge $\mathbf{\phi(\alpha)}$};
 \end{tikzpicture}}\hfill $(ii)$
\end{center}
\begin{center}\scalebox{0.55}{\begin{tikzpicture}[ultra thick]
 \draw[black, thick] (0,0)--(7,0); \draw[black, thick] (0,0)--(0,2);
  \draw[black, line width=3pt] (3,0.4)--(3.4,0);
  \draw[black, dashed, very thick] (0.4,0.4)--(3,0.4);
  \draw[black, line width=3pt] (0,0)--(0.4,0.4);
 \draw[orange,very thick] (0.4,0.4) parabola bend (1.7,1.9) (3,0.4);
  \draw[orange,very thick] (3.4,0) parabola bend (5.1,1.6) (6.8,0);
 \draw  (1.7,1) node {\huge $\mathbf{\alpha}$};
 \draw  (5.1,0.8) node {\huge $\mathbf{\beta}$};

 \end{tikzpicture}}
 \qquad$\stackrel{\phi}{\longrightarrow}$\qquad
\scalebox{0.55}{\begin{tikzpicture}[ultra thick]
 \draw[black, thick] (0,0)--(7,0); \draw[black, thick] (0,0)--(0,2);
  \draw[black, line width=3pt] (0,0)--(0.4,0.4)--(0.8,0.8)--(1.2,0.4);
  \draw[orange,very thick] (1.2,0.4) parabola bend (2.3,1.9) (3.4,0.4);
  \draw[black, dashed, very thick] (0.4,0.4)--(3.4,0.4);
  \draw[black, line width=3pt] (3.4,0.4)--(3.8,0);

   \draw[orange,very thick] (3.8,0) parabola bend (5.3,1.6) (6.8,0);

 \draw[very thick]  (2.3,1) node {\huge $\mathbf{\phi(\alpha)}$};
 \draw  (5.3,0.6) node {\huge $\mathbf{\phi(\beta)}$};
 \end{tikzpicture}}\hfill $(iii)$
\end{center}

\begin{center}\scalebox{0.55}{\begin{tikzpicture}[ultra thick]
 \draw[black, thick] (0,0)--(7,0); \draw[black, thick] (0,0)--(0,2);
  \draw[black, line width=3pt] (3,0.8)--(3.4,0);
  \draw[black, dashed, very thick] (0.4,0.4)--(3,0.8);
  \draw[black, line width=3pt] (0,0)--(0.4,0.4);
 \draw[orange,very thick] (0.4,0.4) parabola bend (1.7,1.9) (3,0.8);
  \draw[orange,very thick] (3.4,0) parabola bend (5.1,1.6) (6.8,0);
 \draw  (1.7,1) node {\huge $\mathbf{\alpha}$};
 \draw  (5.1,0.8) node {\huge $\mathbf{\beta}$};
 \draw (3.5,1.2) node {$k=2$};
 \end{tikzpicture}}
 \qquad$\stackrel{\phi}{\longrightarrow}$\qquad
\scalebox{0.55}{\begin{tikzpicture}[ultra thick]
 \draw[black, thick] (0,0)--(7,0); \draw[black, thick] (0,0)--(0,2);

  \draw[black, line width=3pt] (0,0)--(0.4,0.4);
  \draw[orange,very thick] (0.4,0.4) parabola bend (1.9,1.9) (3.4,0.4);
  \draw[black, dashed, very thick] (0.4,0.4)--(3.4,0.4);
  \draw[black, line width=3pt] (3.4,0.4)--(3.8,0);

   \draw[orange,very thick] (3.8,0) parabola bend (5.3,1.6) (6.8,0);

 \draw[very thick]  (1.9,1) node {\huge $\phi(\alpha D)$};
 \draw  (5.3,0.6) node {\huge $\mathbf{\phi(\beta)}$};

 \end{tikzpicture}}\hfill $(iv)$
\end{center}
\begin{center}\scalebox{0.55}{\begin{tikzpicture}[ultra thick]
 \draw[black, thick] (0,0)--(7,0); \draw[black, thick] (0,0)--(0,2);
  \draw[black, line width=3pt] (3,1.2)--(3.4,0);
  \draw[black, dashed, very thick] (0.4,0.4)--(3,1.2);
  \draw[black, line width=3pt] (0,0)--(0.4,0.4);
 \draw[orange,very thick] (0.4,0.4) parabola bend (1.7,1.9) (3,1.2);
  \draw[orange,very thick] (3.4,0) parabola bend (5.1,1.6) (6.8,0);
 \draw  (1.7,1.2) node {\huge $\mathbf{\alpha}$};
 \draw  (5.1,0.8) node {\huge $\mathbf{\beta}$};
 \draw (3.6,1.4) node {$k\geq 3$};

 \end{tikzpicture}}
 \qquad$\stackrel{\phi}{\longrightarrow}$\qquad
\scalebox{0.55}{\begin{tikzpicture}[ultra thick]
 \draw[black, thick] (0,0)--(7,0); \draw[black, thick] (0,0)--(0,2);

  \draw[black, line width=3pt] (0,0)--(0.4,0.4)--(0.8,0);
  \draw[orange,very thick] (0.8,0) parabola bend (2.3,1.6) (3.8,0);

   \draw[orange,very thick] (3.8,0) parabola bend (5.3,1.6) (6.8,0);

 \draw[very thick]  (2.3,0.5) node {\LARGE $\phi(\alpha D_{k-1})$};
 \draw  (5.3,0.6) node {\huge $\mathbf{\phi(\beta)}$};

 \end{tikzpicture}}\hfill $(v)$
\end{center}
\caption{Illustration of the bijection $\phi$ between $\mathcal{M}_n$ and $\mathcal{A}_n$.}
\label{fig2}\end{figure}
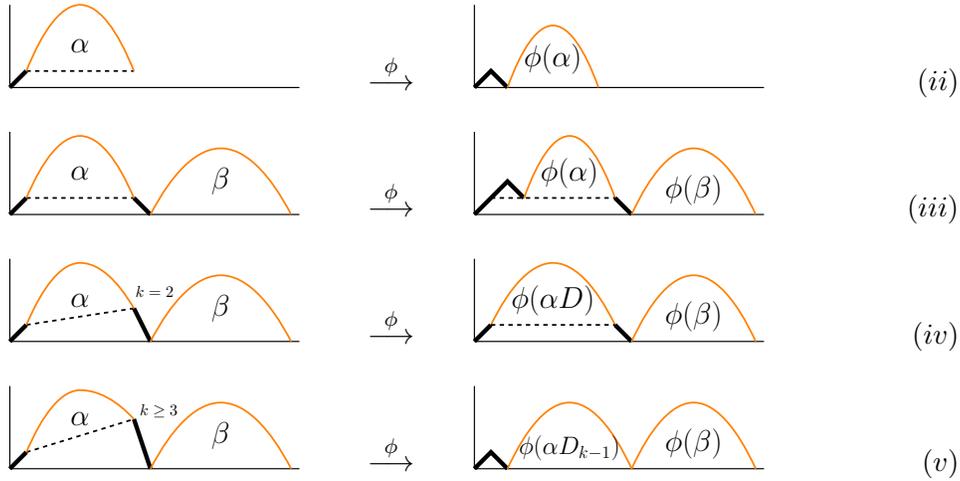

Due to the recursive definition, the image by $\phi$ of a length $n$ Dyck meander with catastrophes
 is a Dyck path of semilength $n$.  For instance, the images of $U$, $UD$, $UUD_2$, $UUDUUD_3$ are respectively $UD$, $UUDD$, $UUUDDD$, $UUDDUUDUUDDD$.
 We refer to Figure \ref{fig2} for an illustration of this mapping.

 \begin{lem} For any $n\geq 0$,

 - if $P\in \mathcal{M}_n$ then we have $\phi(P)\in\mathcal{A}_n$,

 - if $P\in \mathcal{D}_{n}$ then we have $\phi(P)\in\mathcal{A}_{2n}^*$, where $\mathcal{A}_{0}^*=\{\epsilon\}$ and for $n\geq 1$, the set $\mathcal{A}_{2n}^*$ consists of semilength $2n$ Dyck paths avoiding the patterns $UUU$ and $DUD$ and starting with $UUD$.
 \label{lem1}
 \end{lem}
 \noindent {\it Proof.} We proceed by induction on $n$. The case $n=0$ is obvious. For $k\leq n$, we assume that for any $P\in \mathcal{D}_k$ we have $\phi(P)\in \mathcal{A}_{2k}^*$ and for any $P\in\mathcal{M}_{k}$ we have $\phi(P)\in\mathcal{A}_{k}$. Now, let us  prove the result  for $k=n+1$.

Whenever $P\in \mathcal{D}_{n+1}$ we can write $P=U\alpha D\beta$ where $\alpha,\beta\in\mathcal{D}$. Thus, we have $\phi(P)=UUD\phi(\alpha)D\phi(\beta)$, and using the recurrence hypothesis on $\alpha$ and $\beta$, $\phi(P)$ is of semilength $2n+2$, starts with $UUD$  and  avoids the pattern $UUU$. Moreover $\phi(\alpha)$ (resp. $\phi(\beta)$) is either empty or it starts with $UUD$, which implies that $\phi(P)$ avoids $DUD$, and thus $\phi(P)\in\mathcal{A}_{2n+2}^*$.

Now let us assume $P\in\mathcal{M}_{n+1}$.

- If $P=U\alpha$ with $\alpha\in \mathcal{M}_n$, then $\phi(P)=UD\phi(\alpha)$ and the recurrence hypothesis implies that $\phi(P)$ avoids $UUU$ and $DUD$ at height $h>0$.

- If $P=U\alpha D\beta$ where $\alpha\in\mathcal{D}$  and $\beta\in\mathcal{M}$, then the first part of the proof implies that $\phi(\alpha)$ avoids $UUU$ and $DUD$, and with the recurrence hypothesis on $\beta$, $\phi(P)=UUD\phi(\alpha)D\phi(\beta)$ belongs to $\mathcal{A}$.

- If $P=U\alpha D_2\beta$ and $\alpha D\in\mathcal{D}$, then using the first part of the proof $\phi(\alpha D)$ is not empty and avoids $UUU$ and $DUD$. The recurrence hypothesis implies that $\phi(P)=U\phi(\alpha D)D\phi(\beta)$ belongs to $\mathcal{A}$.

- If $P=U\alpha D_i\beta$ where $i\geq 3$ and $\alpha D_{i-1}$ ends on the $x$-axis, then using a simple  induction on $i\geq 2$, $\phi(\alpha D_{i-1})$ is not empty and avoids $UUU$ and $DUD$. The recurrence hypothesis implies that $\phi(P)=U\phi(\alpha D_{i-1})D\phi(\beta)\in\mathcal{A}$.

The induction is completed.
 \hfill $\Box$

\begin{thm} For $n\geq 0$, the map $\phi:\mathcal{M}_n\rightarrow\mathcal{A}_n$  is a bijection. Moreover, we have $\phi(\mathcal{D}_{n})=\mathcal{A}_{2n}^*$.
\label{thm1}
\end{thm}

\noindent {\it Proof.}  Due to the enumerative results in \cite{Ban} (see  Corollary 2.4) and the above lemma,  it suffices to prove that $\phi$ is injective  from $\mathcal{M}_{n}$ to $\mathcal{A}_{n}$.
We proceed by induction on $n$. The case $n=0$ is obvious. For $k\leq n$, we assume that $\phi$ is an injection from $\mathcal{M}_k$ to $\mathcal{A}_k$, and we prove the result  for $k=n+1$.

According to the  definition of $\phi$ and Lemma~\ref{lem1}, the image by $\phi$ of  $P\in\mathcal{M}$  satisfying
($ii$) is a Dyck path starting by $(UD)^kR$ for some $k\geq 1$ where $R$ is a Dyck path in $\mathcal{A}^*_{2i}$ for some $i\geq 0$, which means that $R$ avoids $UUU$; a meander satisfying ($iii$) is sent by $\phi$ to a Dyck path in  $\mathcal{A}^*_{2i}$ for some $i\geq 1$; a meander satisfying ($iv$) is sent to a Dyck path starting with $UUUD$; and a meander satisfying ($v$) is sent to a Dyck path starting with $(UD)^k$ for some $k\geq 1$ and such that it contains an occurrence $UUU$ on the $x$-axis.  Then, for $P,Q\in\mathcal{M}_{n+1}$, $\phi(P)=\phi(Q)$ implies that $P$ and $Q$ belong to the same case ($i$), ($ii$), ($iii$), ($iv$) or ($v$).  So, the recurrence hypothesis induces $P=Q$ which completes the induction. Thus $\phi$ is injective.
 Since $\mathcal{M}_n$ and $\mathcal{A}_n$ have the same cardinality (see \cite{Ban} and \href{https://oeis.org/A274115}{A274115} in \cite{Sloa}), $\phi$ is a bijection.

 Considering the previous lemma, it suffices to check that $\mathcal{A}_{2n}^*$ is counted by the Catalan numbers in order  to prove  $\phi(\mathcal{D}_{n})=\mathcal{A}_{2n}^*$. A Dyck path $P\in \mathcal{A}_{2n}^*$ is either empty or  it
 consists of a sequence of $UUD\alpha D$ where  $\alpha$ belongs to $\mathcal{A}^*_{2n-2}$. Let $A^*(x)$ be  the generating function for the cardinality of  $\mathcal{A}_{2n}^*$ (with respect to the semilength).  We obtain the following functional equations $A^*(x)=1+\frac{x^2A^*(x)}{1-x^2A^*(x)}$ which implies that $\mathcal{A}_{2n}^*$ is counted by the $n$th Catalan number. Therefore $\phi : \mathcal{D}_{n}\rightarrow\mathcal{A}_{2n}^*$ is a bijection.
\hfill $\Box$

\begin{rem} In \cite{Mane}, it is proven that the set $\mathcal{A}_n$ is a representative set of the equivalence classes  modulo the pattern $DUU$ on Dyck paths, i.e. two Dyck paths $P$ and $Q$ are equivalent if and only
 if the positions of the occurrences $DUU$ are the same in $P$ and $Q$ (see also \cite{Bar}). So, the bijection $\phi$  establishes a direct correspondence between these classes and Dyck meanders with catastrophes.
\end{rem}
\bigskip
Let $\mathcal{A}'_n$ be the subset of $\mathcal{A}_n$ consisting of paths $P$ such that  any occurrence $UD$ on the $x$-axis in $P$  appears before an occurrence  of $UUU$  (not necessarily contiguous to the occurrence $UD$). The next theorem gives a bijection between $\mathcal{A}'_{n}$ and the set $\mathcal{E}_n$ of length $n$ Dyck excursions with catastrophes.

\begin{thm} For $n\geq 0$, we have $\phi(\mathcal{E}_n)=\mathcal{A}'_{n}$.
\label{thm2}
\end{thm}

\noindent {\it Proof.} Thanks to Theorem \ref{thm1}, it suffices to check that for any $P\in \mathcal{E}_n$, $\phi(P)\in \mathcal{A}'_{n}$, and $|\mathcal{A}'_{n}|=|\mathcal{E}_n|$. Any $P\in \mathcal{E}_n$  satisfies one of the cases $(i)$, $(iii)$, $(iv)$ and $(v)$ with $\beta\in\mathcal{E}$. We proceed by induction on the length in order to prove that $\phi(P)\in \mathcal{A}'_{n}$. The case $(i)$ is obvious. Whenever $P$ satisfies the cases $(iii)$ or $(iv)$, the only possibility for an occurrence of  $UD$ to appear at height zero in  $\phi(P)=UUD\phi(\alpha)D\phi(\beta)$ (resp. $\phi(P)=U\phi(\alpha D)D\phi(\beta)$) is 
 to be inside $\phi(\beta)$. Applying the recurrence hypothesis on $\beta$,  $\phi(P)\in \mathcal{A}'_{n}$.
For a path $P$ satisfying the case $(v)$, we have seen in the proof of Theorem \ref{thm1} that $\phi(P)$ starts necessarily with $(UD)^k$ for $k\geq 1$  followed by $UUU$. Using the recurrence hypothesis for $\beta$, we obtain $\phi(P)\in \mathcal{A}'_{n}$. The induction is completed.

Now, let us prove that $|\mathcal{A}'_{n}|=|\mathcal{E}_n|$. Any path $P\in\mathcal{A}'_n$ satisfies one of the following two cases: ($a$) $P\in\mathcal{A}_{n}$ does not contain any occurrence $UD$ on the $x$-axis, and ($b$) $P=QUDR$ where $Q\in\mathcal{A}$ and  $R\in\mathcal{A}$ such that $R$ contains at least one occurrence of $UUU$ and avoids any occurrence $UD$ on the $x$-axis.
Let $\mathcal{K}$ (resp. $\overline{\mathcal{K}}$) be the set of Dyck paths in $\mathcal{A}$ satisfying ($a$) (resp. ($b$)), and let $K(x)$ and $\overline{K}(x)$ be the corresponding  generating functions for their cardinalities with respect to the semilength. Obviously, the generating function $A'(x)$ for $\mathcal{A}'=\cup_{n\geq 0}\mathcal{A}_n$ satisfies 
$$A'(x)=K(x)+\overline{K}(x).$$

A nonempty path $P\in\mathcal{K}$ can be  decomposed  $P=U\alpha D \beta$ where $\beta\in\mathcal{K}$ and $\alpha$ is a nonempty Dyck path avoiding $UUU$ and $DUD$. Then either $\alpha\in \mathcal{A}^\star\backslash\{\epsilon\}$ or $\alpha=UD\alpha'$  with  $\alpha'\in \mathcal{A}^\star$.
Thus the generating function for $\mathcal{K}$ is given by 
$$K(x)=1+x(A^\star(x)-1+xA^\star(x))\cdot K(x).$$

 Due to the form of a path $P\in\overline{\mathcal{K}}$, we deduce the  functional equation 
$$\overline{K}(x)=A(x)x R(x)$$

where $A(x)$ is the generating function for $\mathcal{A}$ and $R(x)$ is the generating function for the paths in $\mathcal{A}$ avoiding any occurrence $UD$ on the $x$-axis and containing at least one occurrence of $UUU$. 
 Due to Theorem~\ref{thm1}, $A(x)$ is also the o.g.f. for Dyck meanders with
 catastrophes that is $A(x) = M(x) = \frac{2 x}{2 x + \left(x + 1\right) 
 \left(\sqrt{1 - 4 x^{2}} - 1\right)}$.
Then, we have $R(x)=K(x)-L(x)$ where $L(x)$ is the generating function for the set $\mathcal{L}$ of Dyck paths in $\mathcal{A}$ avoiding $UD$ and $UUU$ on the $x$-axis. Note that $\mathcal{L}$ is exactly the set $\mathcal{A}^\star$, then $$R(x)=K(x)-A^\star(x).$$

Combining the previous equations, we obtain $A'(x)={\frac {2-3\,x-2\,{x}^{2}+x
\sqrt {1-4\,{x}^{2}}}
{2-2x - 4x^2 - 2x^3}
}$ which is exactly the generating function of $\mathcal{E}_n$ found by \cite{Ban}.
\hfill $\Box$


\section{Motzkin meanders with catastrophes}

In this section we exhibit a constructive bijection between the set $\mathcal{M}'_n$  of length $n$ Motzkin meanders with catastrophes and the set
 $\mathcal{B}_{n+1}$ of semilength $n+1$ Dyck paths avoiding the patterns $UUU$ at height $h\geq 2$.
Before defining this bijection we recall that there exists a one-to-one correspondence  $\chi$ between length $n$ Motzkin paths and semilength $n$ Dyck paths avoiding $UUU$. From a Dyck path avoiding $UUU$, we replace each $UUD$ with $U$, and we replace each remaining $UD$ with $F$. For instance, the image by $\chi$ of $UFUDD$ is $UUDUDUUDDD$ (see \cite{Cal,Eli}).

Now, let us use $\chi$ in order to define recursively the map $\psi$ from $\mathcal{M}'$ to $\mathcal{D}$ as follows.
For $P\in\mathcal{M}'$, we set $\psi(P)=$
 $$\left\{\begin{array}{ll}
 UD& \mbox{ if } P=\epsilon,\\
 UD\psi(\alpha) &\mbox{ if } P=F\alpha,\\
 UU\chi(\alpha_1 )DU\chi(\alpha_2)D \ldots U\chi(\alpha_k)DD  & \mbox{ if } P= U\alpha_1 U \alpha_2 \ldots U\alpha_k,\\
 UU\chi(\alpha_1 )DU\chi(\alpha_2)D \ldots U\chi(\alpha_k)DD\psi(\beta) & \mbox{ if } P= U\alpha_1 U \alpha_2 \ldots U\alpha_kD_k\beta,\\
\end{array}\right.$$
where $k\geq 1$, $\alpha,\beta\in\mathcal{M}'$, and $\alpha_1,\alpha_2,\ldots , \alpha_k$ are some possibly empty Motzkin  paths (considering that $D_1$ is defined to be $D$).

\begin{figure}[h]
\begin{center}\scalebox{0.45}{\begin{tikzpicture}[ultra thick]
 \draw[black, thick] (0,0)--(9,0); \draw[black, thick] (0,0)--(0,2);
  \draw[black, line width=3pt] (0,0)--(0.4,0);
 \draw[orange,very thick] (0.4,0) parabola bend (1.7,1) (3,0);
 \draw  (1.7,0.5) node {\huge $\alpha$};
 \end{tikzpicture}}
 \qquad$\stackrel{\psi}{\longrightarrow}$\qquad
\scalebox{0.45}{\begin{tikzpicture}[ultra thick]
 \draw[black, thick] (0,0)--(11.1,0); \draw[black, thick] (0,0)--(0,2);
 \draw[black, line width=3pt] (0,0)--(0.4,0.4)--(0.8,0);
 \draw[orange,very thick] (0.8,0) parabola bend (1.9,1) (3,0);
 \draw  (1.9,0.5) node {\huge $\psi(\alpha)$};
 \end{tikzpicture}}\hfill $(ii)$
\end{center}
\begin{center}\scalebox{0.45}{\begin{tikzpicture}[ultra thick]
 \draw[black, thick] (0,0)--(9,0); \draw[black, thick] (0,0)--(0,2);
 \draw[black, dashed, very thick] (0.4,0.4)--(2,0.4);
  \draw[black, line width=3pt] (0,0)--(0.4,0.4);
 \draw[orange,very thick] (0.4,0.4) parabola bend (1.2,1) (2,0.4);
  \draw[black, dashed, very thick] (2.4,0.8)--(4,0.8);
  \draw[black, line width=3pt] (2,0.4)--(2.4,0.8);
 \draw[orange,very thick] (2.4,0.8) parabola bend (3.2,1.4) (4,0.8);
 \draw[black, dashed, very thick] (4,0.8)--(4.8,1.6);
 \draw[black, line width=3pt] (4.8,1.6)--(5.2,2);
 \draw[black, dashed, very thick] (5.2,2)--(6.8,2);
 \draw[orange,very thick] (5.2,2) parabola bend (6,2.6) (6.8,2);

 \draw  (1.2,0.6) node {\huge $\alpha_1$};
  \draw  (3.2,1) node {\huge $\alpha_2$};
   \draw  (6,2.2) node {\huge $\alpha_k$};
 \end{tikzpicture}}
 \qquad$\stackrel{\psi}{\longrightarrow}$\qquad
\scalebox{0.45}{\begin{tikzpicture}[ultra thick]
 \draw[black, thick] (0,0)--(11,0); \draw[black, thick] (0,0)--(0,2);
 \draw[black, line width=3pt] (0,0)--(0.4,0.4)--(0.8,0.8);
 \draw[black, dashed, very thick] (0.8,0.8)--(2.4,0.8);

 \draw[orange,very thick] (0.8,0.8) parabola bend (1.6,1.4) (2.4,0.8);
  \draw[black, line width=3pt] (2.4,0.8)--(2.8,0.4)--(3.2,0.8);
 \draw[orange,very thick] (3.2,0.8) parabola bend (4,1.4) (4.8,0.8);
  \draw[black, line width=3pt] (4.8,0.8)--(5.2,0.4);
  \draw[black, dashed, very thick] (3.2,0.8)--(4.8,0.8);
   \draw[black, dashed, very thick] (5.5,0.4)--(6,0.4);
  \draw[black, line width=3pt] (6.6,0.4)--(7,0.8);
  \draw[orange,very thick] (7,0.8) parabola bend (7.8,1.4) (8.6,0.8);
   \draw[black, dashed, very thick] (7,0.8)--(8.6,0.8);
   \draw[black, line width=3pt] (8.6,0.8)--(9.4,0);
 \draw  (1,1.7) node {\huge $\chi(\alpha_1)$};
 \draw  (3.5,1.7) node {\huge $\chi(\alpha_2)$};
 \draw  (7.3,1.7) node {\huge $\chi(\alpha_k)$};
 \end{tikzpicture}}\hfill $(iii)$
\end{center}
\begin{center}\scalebox{0.45}{\begin{tikzpicture}[ultra thick]
 \draw[black, thick] (0,0)--(9,0); \draw[black, thick] (0,0)--(0,2);
 \draw[black, dashed, very thick] (0.4,0.4)--(2,0.4);
  \draw[black, line width=3pt] (0,0)--(0.4,0.4);
 \draw[orange,very thick] (0.4,0.4) parabola bend (1.2,1) (2,0.4);
  \draw[black, dashed, very thick] (2.4,0.8)--(4,0.8);
  \draw[black, line width=3pt] (2,0.4)--(2.4,0.8);
 \draw[orange,very thick] (2.4,0.8) parabola bend (3.2,1.4) (4,0.8);
 \draw[black, dashed, very thick] (4,0.8)--(4.8,1.6);
 \draw[black, line width=3pt] (4.8,1.6)--(5.2,2);
 \draw[black, dashed, very thick] (5.2,2)--(6.8,2);
 \draw[orange,very thick] (5.2,2) parabola bend (6,2.6) (6.8,2);
\draw[black, line width=3pt] (6.8,2)--(7.2,0);
\draw[orange,very thick] (7.2,0) parabola bend (8,0.6) (8.8,0);

  \draw  (1.2,0.6) node {\huge $\alpha_1$};
  \draw  (3.2,1) node {\huge $\alpha_2$};
   \draw  (6,2.2) node {\huge $\alpha_k$};
 \draw  (8.2,1) node {\huge $\beta$};
 \end{tikzpicture}}
 \qquad$\stackrel{\psi}{\longrightarrow}$\qquad
\scalebox{0.45}{\begin{tikzpicture}[ultra thick]
 \draw[black, thick] (0,0)--(11.1,0); \draw[black, thick] (0,0)--(0,2);
 \draw[black, line width=3pt] (0,0)--(0.4,0.4)--(0.8,0.8);
 \draw[black, dashed, very thick] (0.8,0.8)--(2.4,0.8);

 \draw[orange,very thick] (0.8,0.8) parabola bend (1.6,1.4) (2.4,0.8);
  \draw[black, line width=3pt] (2.4,0.8)--(2.8,0.4)--(3.2,0.8);
 \draw[orange,very thick] (3.2,0.8) parabola bend (4,1.4) (4.8,0.8);
  \draw[black, line width=3pt] (4.8,0.8)--(5.2,0.4);
  \draw[black, dashed, very thick] (3.2,0.8)--(4.8,0.8);
   \draw[black, dashed, very thick] (5.5,0.4)--(6,0.4);
  \draw[black, line width=3pt] (6.6,0.4)--(7,0.8);
  \draw[orange,very thick] (7,0.8) parabola bend (7.8,1.4) (8.6,0.8);
   \draw[black, dashed, very thick] (7,0.8)--(8.6,0.8);
   \draw[black, line width=3pt] (8.6,0.8)--(9.4,0);
   \draw[orange,very thick] (9.4,0) parabola bend (10.2,0.6) (11,0);

 \draw  (1,1.7) node {\huge $\chi(\alpha_1)$};
 \draw  (3.5,1.7) node {\huge $\chi(\alpha_2)$};
 \draw  (7.3,1.7) node {\huge $\chi(\alpha_k)$};

 \draw  (10.2,1) node {\huge $\mathbf{\psi(\beta)}$};
 \end{tikzpicture}}\hfill $(iv)$
\end{center}
\caption{Illustration of the bijection $\psi$ between $\mathcal{M}'_n$
and $\mathcal{B}_{n+1}$.
}
\label{fig3}
\end{figure}
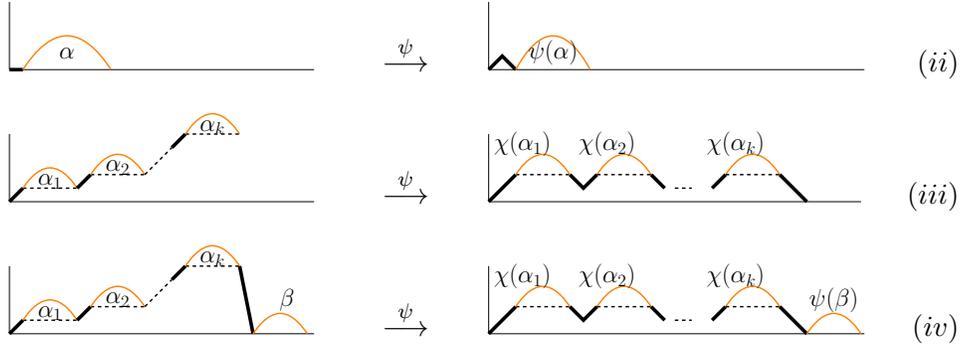

Clearly, the image by $\psi$ of a length $n$ Motzkin meander with catastrophes
 is a Dyck path of semilength $n+1$.  For instance, the images by $\psi$ of $\epsilon$, $F$, $UD$, $UUD_2$, $UUDUUD_3$ are respectively $UD$, $UDUD$, $UUDDUD$, $UUDUDDUD$,
 $UUUUDDUDUDDUD$. We refer to Figure \ref{fig3} for an illustration of this mapping.

A simple observation provides the following results.

\begin{thm} For $n\geq 0$, the map $\psi$, defined above, induces a bijection from $\mathcal{M}'_n$ to $\mathcal{B}_{n+1}$. Moreover, the image of the set of length $n$ Motzkin paths is the set of semilength $n+1$ Dyck paths avoiding $UUU$ at height $h\geq 2$ and the pattern $DU$ at  height one.\end{thm}

\begin{cor} For $n\geq 0$, $\psi(\mathcal{E}'_n)$ is the set of Dyck paths in $\mathcal{B}_{n+1}$ ending with $UD$,  which implies that $\psi$  induces a one-to-one correspondence from paths $P\in\mathcal{E}'_n$  to $\mathcal{B}_{n}$ after deleting the  last two steps $UD$ from $\psi(P)$.
\end{cor}

\section{Acknowledgement} We would like to thank Cyril Banderier for suggesting us to explore constructive bijections between  meanders with catastrophes and  Dyck paths avoiding some patterns.

\end{document}